\documentstyle[12pt,a4]{article}

\title{
  \begin{flushright}
    \normalsize q-alg/9506016\\
    RIMS-1019\\
    YITP/K-1111\\
    June 16, 1995      
  \end{flushright}\medskip
  Level-$0$ structure of level-$1$ $\uqlarge$-modules\\
  and Macdonald polynomials}

\author{
Michio Jimbo\thanks{Department of Mathematics, Faculty of Science,
                            Kyoto University, Kyoto 606, Japan.}\quad
Rinat Kedem\thanks{Research Institute for Mathematical Sciences,
                            Kyoto University, Kyoto 606, Japan.}\quad
Hitoshi Konno\thanks{Yukawa Institute for Theoretical Physics,
                            Kyoto University, Kyoto 606, Japan.}\\
Tetsuji Miwa$^\dagger$\quad Jens-Ulrik H. Petersen$^\dagger$}

\date{}

\newfont{\twlvmsb}{msbm10 scaled\magstep1}
\newfont{\ninemsb}{msbm9}
\newfont{\sixmsb}{msbm6}
\newfam\msbfam
\textfont\msbfam=\twlvmsb
\scriptfont\msbfam=\ninemsb
\scriptscriptfont\msbfam=\sixmsb
\catcode`\@=11
\def\Bbb{\ifmmode\let\next\Bbb@\else
  \def\next{\errmessage{Use \string\Bbb\space only in math mode}}\fi\next}
\def\Bbb@#1{{\Bbb@@{#1}}}
\def\Bbb@@#1{\fam\msbfam#1}
\newfont{\largeeufm}{eufm10 scaled\magstep4}
\newfont{\twlveufm}{eufm10 scaled\magstep1}
\newfont{\elveufm}{eufm10 at 11pt}
\newfont{\teneufm}{eufm10}
\newfont{\nineeufm}{eufm9}
\newfam\eufam
\textfont\eufam=\twlveufm
\scriptfont\eufam=\nineeufm
\def\frak{\ifmmode\let\next\frak@\else
  \def\next{\errmessage{Use \string\frak\space only in math mode}}\fi\next}
\def\frak@#1{{\fam\eufam{{#1}}}}
\catcode`@=12
\newcommand{\uqlarge}{U_q(\widehat{\mbox{\largeeufm sl}}_2)}
\newcommand{\uqsmall}{U_q(\widehat{\mbox{\elveufm sl}}_2)}

\def\be{\begin{equation}}
\def\en{\end{equation}}
\def\bea{\begin{eqnarray}}
\def\ena{\end{eqnarray}}
\def\bean{\begin{eqnarray*}}
\def\enan{\end{eqnarray*}}
\def\lb#1{\label{eqn:#1}}
\def\rf#1{(\ref{eqn:#1})}
\def\qed{\hfill\fbox{}}
\def\pf{\noindent{\it Proof.\quad}}
\def\displ#1{{\displaystyle #1}}

\def\mod{{\rm mod}\,}
\newcommand{\gsl}{{\frak{sl}}}
\newcommand{\ggl}{{\frak{gl}}}
\newcommand{\slt}{\gsl_2}
\newcommand{\slh}{\widehat{\gsl}}
\newcommand{\glh}{\widehat{\ggl}}
\def\slth{\slh_2}
\def\slnh{\slh_n}
\def\glth{\glh_2}

\def\op{\scriptstyle {\rm op}}
\def\V{{\cal V}}
\def\uq{U_q(\slth)}
\def\uqc{U_q(\slt)}

\def\End{{\rm End}}
\def\Ker{{\rm Ker}}

\def\ch{{\rm ch}\,}
\def\tr{{\rm tr}\,}
\def\ad{{\rm ad}\,}

\def\ket#1{|#1\rangle}

\def\vaf{V_{\hbox{\scriptsize aff}}}

\def\H{{\cal H}}
\def\Z{{\Bbb Z}}
\def\Q{{\Bbb Q}}
\def\C{{\Bbb C}}
\def\e{\varepsilon}
\def\vep{\varepsilon}

\def\cL{{\cal L}}
\def\tL{\tilde{\cal L}}
\def\bL{\overline{L}}

\def\bR{\overline{R}}

\def\cN{{\cal N}}
\def\tN{\tilde{\cN}}
\def\hV{\hat{\cal V}}
\def\Vh{\hat{\cal V}}
\newcommand{\rh}{\hat{\rho}}

\newtheorem{prop}{Proposition}
\newtheorem{thm}[prop]{Theorem}
\newtheorem{lem}[prop]{Lemma}

\begin{document}
\maketitle
\begin{abstract}
  The level-$1$ integrable highest weight modules of
  $\uqsmall$ admit a level-$0$ action of the same algebra.
  This action is defined using the affine Hecke algebra and the basis
  of the level-$1$ module generated by components of vertex operators.
  Each level-$1$ module is a direct sum of finite-dimensional
  irreducible level-$0$ modules, whose highest weight vector is
  expressed in terms of Macdonald polynomials.  This decomposition
  leads to the fermionic character formula for the level-$1$ modules.
\end{abstract}

\setcounter{equation}{0}
\setcounter{prop}{0}
\section{Introduction}\label{sec:intro}

The present article is concerned with a curious finding in
representation theory, noted recently in connection with conformal
field theory (CFT) and solvable lattice models.  In the context of
quantum affine algebras, the statement is the following.  Let
$V(\Lambda_i)$ ($i=0,1$) be a level-$1$ integrable highest weight
module of $\uq$.  Then it admits a level-$0$ action of the same
algebra $\uq$.  This second action leaves invariant each homogeneous
component $V(\Lambda_i)_{-n}$, and the whole module becomes a direct
sum of finite-dimensional irreducible constituents, which can be
described explicitly.

The origin of this observation goes back to the study of spin chains
with long-range interaction (the Haldane-Shastry model)~\cite{HHTBP}.
This model has a remarkable property that the Yangian $Y(\gsl_n)$ acts
as an exact symmetry even for finite chains.  By considering the
continuum limit, it was suggested in~\cite{HHTBP}, and subsequently
confirmed in~\cite{BGHP,BPS,Schou,BLSa}, that there is an action of
the Yangian on level-$1$ integrable modules of the affine Lie algebra
$\slnh$.  This action is related to the fermionic character formula
for the level-$1$ module conjectured in~\cite{Mel}.  Similar results
are expected to hold for higher level representations as
well~\cite{BLSb}.  The fermionic expressions conjectured there for the
characters of the level-$k$ modules, subsequently proved
in~\cite{NaYa,ANOT}, strongly support the validity of this picture.

What we discuss in this paper is a $q$-deformation of these structures
in the simplest case of $\uq$ with level-$1$.  Namely we consider
level-$1$ modules of $\uq$, and define a level-$0$ action on them.
The quantum affine algebra with level $0$, plays the role of the
Yangian in CFT.  It is not yet clear whether such a $q$-deformation is
related to some physical models like the Haldane-Shastry chain.  Apart
from the technical complexity, the essence of the construction does
not differ very much from the conformal case.  In this respect we are
not claiming any sort of methodical novelty.  Our aim here is to
supply the mathematical details and give a coherent account of this
yet mysterious phenomenon.

Here is the plan of the paper.  In Section~\ref{sec:action} we
introduce a level-$0$ action of $\uq$ on the sum of level-$1$
integrable highest weight modules $\H=V(\Lambda_0)\oplus
V(\Lambda_1)$.  This part is a review of our previous
paper~\cite{USC}.  The level-$0$ action is defined with the help of
the affine Hecke algebras which appear in~\cite{BGHP}, and the action
is given on certain generating series of the vectors in $\H$.  To
ensure the construction is well-defined, we introduce appropriate
completions of the spaces, and sort out all the linear relations among
the generating series mentioned above.  In the conformal case, the
Yangian generators have explicit realizations in terms of currents.
Lacking such formulas in the deformed case, our construction is rather
indirect.  In Section~\ref{sec:hwv} we define a family of highest
weight vectors $\omega_{\lambda,N}\in \H$ with the aid of Macdonald
polynomials, following the ideas in~\cite{BPS}.  These vectors are
indexed by $N\in\Z_{\ge 0}$ and a partition $\lambda$ of length at
most $N$.  We remark that the vectors $\omega_{\lambda,N}$ can be
given an explicit formula in terms of bosons and Macdonald symmetric
functions (see~\rf{hwvB}).  This point will not be used in the rest of
the paper.  We show in Section~\ref{sec:decomp} that
$\omega_{\lambda,N}$'s generate inside $\H$ irreducible modules of the
form (see~(\ref{eqn:kn}, \ref{eqn:str}, \ref{eqn:Wln}))
\[
W_{\lambda,N}=W_{n_1}(a_1)\otimes \cdots \otimes W_{n_m}(a_m).
\]

Here $W_n(a)$ is the $(n+1)$-dimensional evaluation module (see
Appendix~\ref{app:Drinf} for conventions).  Moreover $\H$ is a direct
sum of them:
\[
V(\Lambda_i)=
\bigoplus_{N\equiv i\,(\mod 2) \atop l(\lambda)\le N}W_{\lambda,N}.
\]
This decomposition corresponds to the known fermionic sum formula for
the characters
\[
\ch V(\Lambda_i)
=\hbox{tr}_{V(\Lambda_i)}\left(q^d\right)
=\sum_{N\equiv i\,(\mod 2) \atop l(\lambda)\le N}
q^{(N^2-i)/4+|\lambda|}\prod_{j=1}^m(n_j+1),
\]
which we use in the course of the proof of these statements.

Throughout this paper we fix a complex number $q$ such that $0<|q|<1$.
In Section~\ref{sec:decomp} we assume that $q$ is transcendental over
$\Q$.

\subsection{Acknowledgements}
We thank Katsuhisa Mimachi, Tomoki Nakanishi, Atsushi Nakayashiki,
Akihiro Tsuchiya and Yasuhiko Yamada for discussions.  This work is
partly supported by Grant-in-Aid for Scientific Research on Priority
Areas 231, the Ministry of Education, Science and Culture.  H.~K. is
supported by Soryushi Shyogakukai.  R.~K. and J.-U.~H.~P. are
supported by the Japan Society for the Promotion of Science.

\setcounter{equation}{0}
\setcounter{prop}{0}
\def\hY{\hat Y}

\section{Level-$0$ Action} \label{sec:action}
In this section
we define the level-$0$ action on the level-$1$
integrable modules of $\uq$.
The proofs of the statements in this section are given in~\cite{USC}.

\subsection{Spinon Basis}
This subsection contains some preliminary material.
Let us fix the notation following~\cite{JM}.
The quantum affine algebra $U=\uq$ is generated by
$e_i,f_i,t_i=q^{h_i}$ ($i=0,1$) and $q^d$.
When we refer to the subalgebra $\uqc$, we mean
the one generated by $e_1,f_1$ and $t_1$.
By $U'$ we mean the subalgebra of $U$ generated by
$e_i$, $f_i$, $t_i$ $(i=0,1)$.
We use the coproduct $\Delta$ given by
\begin{eqnarray*}
&&\Delta(e_i)=e_i\otimes 1+t_i\otimes e_i,
\qquad
\Delta(f_i)=f_i\otimes t_i^{-1}+ 1\otimes f_i,
\\
&&\Delta(q^h)=q^h\otimes q^h
\qquad (h=h_i,d).
\end{eqnarray*}

Let $V=\C v_+\oplus \C v_-$ be the irreducible two-dimensional
module of $U_q(\slt)$.
Let further
\[
\vaf={\rm span}_\C \{ v_{\vep,n}\mid \vep=\pm,n\in \Z\,\}
\]
be the level-$0$ $U$-module defined as follows.
\begin{eqnarray*}
&&e_0v_{+,n}=v_{-,n+1},~e_0v_{-,n}=0,
\quad f_0v_{+,n}=0,~f_0 v_{-,n}=v_{+,n-1},\\
&&e_1v_{+,n}=0,~e_1 v_{-,n}=v_{+,n},
\quad f_1 v_{+,n}=v_{-,n},~f_1v_{-,n}=0,\\
&&t_0^{-1} v_{\pm,n}=t_1 v_{\pm,n}=q^{\pm 1}v_{\pm,n},
\quad
q^d v_{\pm,n}=q^n v_{\pm,n}.
\end{eqnarray*}
We shall use the generating series
\[
v_{\vep}(z)=\sum_{n}v_{\vep,n}z^{-n}.
\]

Let $V(\Lambda_i)=U\ket{i}$ $(i=0,1)$ be the integrable
level-$1$ $U$-module with highest weight vector $\ket{i}$.
The index $i$ is to be read modulo $2$, i.e., $\Lambda_{i+2}=\Lambda_i$.
Set $\H=V(\Lambda_0)\oplus V(\Lambda_1)$.
This is a $\Z$-graded vector space
\[
\H=\bigoplus_{r\le 0}\H_r,
\qquad
\H_r=\{u\in\H|q^du=q^ru\}.
\]
In particular, $\H_0=\C\ket{0}\oplus\C\ket{1}$.
We regard $\End_\C\H$ as a level-$0$ $U$-module
by the opposite adjoint action
\[
\ad\!^{\op} x.f=\sum x_{(2)}~f~ a^{-1}\bigl(x_{(1)}\bigr),
\]
where $f\in\End_\C\H$, $x\in U$, $\Delta(x)=\sum x_{(1)}\otimes x_{(2)}$
and $a$ denotes the antipode.

The level-$0$ action will be defined with the aid of
the type-I vertex operator
\[
\tilde{\Phi}^*_\vep(z)=\sum_{n\in\Z}\tilde{\Phi}^*_{\vep,n}z^{-n}.
\]
By definition, it is the generating series of operators
$\tilde{\Phi}^*_{\vep,n}\in\End_\C\H$ with the following properties:
(i)~The map $u\otimes v_{\vep,n}\mapsto\tilde{\Phi}^*_{\vep,n}u$ from
$\H\otimes\vaf$ to $\H$ is an intertwiner of $U$-modules,
(ii)~$\tilde{\Phi}^*_{\vep,n}\H_r\subset\H_{r+n}$,
(iii)~$\tilde{\Phi}^*_{+,0}\ket{0}=\ket{1}$,
$\tilde{\Phi}^*_{-,0}\ket{1}=\ket{0}$.

Consider the vector space
\[
\V=\bigoplus_{N\ge 0}\V_N,\qquad
\V_N=\vaf^{\otimes N},
\]
with the $U$-action given by the {\it opposite} coproduct
$\Delta^{\scriptstyle {\rm op}}(x)=\sum x_{(2)}\otimes x_{(1)}$.

\begin{prop}
The map
\begin{eqnarray}
&&\rho_0~:~\V~~ \longrightarrow ~~\H,
\label{eqn:rhoz}\\
&&v_{\vep_1,n_1}\otimes \cdots \otimes v_{\vep_N,n_N}
{}~~\mapsto~~
\tilde{\Phi}^*_{\vep_1,n_1}\cdots \tilde{\Phi}^*_{\vep_N,n_N}\ket{0},
\label{eqn:spinon}
\end{eqnarray}
is surjective and $\uqc$-linear.
\end{prop}
The vectors of the form~\rf{spinon} are not linearly independent.  As
shown in~\cite{Notes}, it is possible to choose a suitable subset to
construct a basis.  A similar basis appeared also in CFT under the
name of the spinon basis, which was used for constructing the Yangian
action on conformal blocks~\cite{BPS,Schou,BLSb}.  We follow the same
strategy in the $q$-deformed case in order to define the level-$0$
action.  For this purpose we need to determine all the linear
relations among~\rf{spinon}~\cite{Notes}.  This will be discussed in
the next subsection.

\subsection{Linear Relations}
The map $\rho_0$ factors through the map
\begin{eqnarray}
&& \rho~:~\V~~ \longrightarrow ~~\End_\C\H
\label{eqn:rho}\\
&&v_{\vep_1,n_1}\otimes \cdots \otimes v_{\vep_N,n_N}
{}~~\mapsto~~
\tilde{\Phi}^*_{\vep_1,n_1}\cdots \tilde{\Phi}^*_{\vep_N,n_N}.
\end{eqnarray}
Among the linear relations for the vectors~\rf{spinon},
the major ones derive from the commutation relations of vertex operators.
The latter involve infinite sums, and
to justify them we need to complete the space $\V$.

Set
\[
\V^{(r)}_N={\rm span}_\C
\{ v_{\vep_1,m_1}\otimes\cdots\otimes v_{\vep_N,m_N}
\mid m_1+\cdots+m_N=r\,\}.
\]
Define a filtration of $\V^{(r)}_N$ by
\begin{eqnarray*}
&&\V^{(r)}_N\supset \cdots \supset \V^{(r)}_N[l]
\supset \V^{(r)}_N[l+1] \supset \cdots,\\
&&\V^{(r)}_N[l]={\rm span}_\C
\{v_{\vep_1,m_1}\otimes\cdots\otimes v_{\vep_N,m_N}\in \V^{(r)}_N\\
&&\quad\mid \max(m_1+m_2+\cdots+m_N,\ldots,m_{N-1}+m_N,m_N)\ge l\,\}.
\end{eqnarray*}
We complete the space $\V^{(r)}_N$ with respect to this filtration
\[
\Vh^{'(r)}_N=\lim_{\longleftarrow\atop l}
\V^{(r)}_N/\V^{(r)}_N[l],
\]
and set $\Vh'_N=\oplus_{r\in\Z}\Vh^{'(r)}_N$,
$\Vh'=\oplus_{N\geq 0}\Vh'_N$.
Since $\H_t=0$ for $t>0$, we have
\[
\rho\left(\V^{(r)}_N[l]\right)\H_s=0\quad \hbox{ if $l+s>0$}.
\]
Therefore, the map $\rho$ and $\rho_0$ extend to the completion:
\[
\rh'~:~\Vh'\longrightarrow\End_\C\H,
\quad
\rh'_0~:~\Vh'\longrightarrow \H.
\]

To state the commutation relations, we
prepare the $R$-matrix $\tilde R(z)\in\End_\C V\otimes V$
\begin{eqnarray}
&&\tilde R(z)v_\vep\otimes v_\vep={z-q^2\over1-q^2z}v_\vep\otimes v_\vep,
\nonumber\\
&&\tilde R(z)v_+\otimes v_-
={(1-q^2)z\over1-q^2z}v_+\otimes v_-+{q(z-1)\over1-q^2z}v_-\otimes v_+,
\nonumber\\
&&\tilde R(z)v_-\otimes v_+
={q(z-1)\over1-q^2z}v_+\otimes v_-+{1-q^2\over1-q^2z}v_-\otimes v_+.
\lb{Rtilde}
\end{eqnarray}
These formulas are regarded as a power series in $z$.
Let $\pi_j:\End_\C V\rightarrow\End_\C V^{\otimes N}$ be the natural injection
to the $j$-th component of $\End_\C V^{\otimes N}\simeq
\End_\C V\otimes\cdots\otimes\End_\C V$. We set
$\tilde R_{j,k}(z)=(\pi_j\otimes\pi_k)\Bigl(\tilde R(z)\Bigr)$.

Consider the generating series in $\Vh'_N$,
\begin{eqnarray}
F_{\vep_1,\ldots,\vep_N}(z_1,\ldots,z_N)
&=&\frac{\prod_{j=1}^N z_j^{(N-j-p_j+p_N)/2}}{\prod_{j<k}\eta(z_k/z_j)}
v_{\vep_1}(z_1)\otimes\cdots\otimes v_{\vep_N}(z_N)
\nonumber\\
&=&
\sum_{m_1,\ldots,m_N\in\Z}
F_{\vep_1\cdots\vep_N,m_1\cdots m_N}
z_1^{-m_1}\cdots z_N^{-m_N},
\label{eqn:gen}
\end{eqnarray}
where
\[
\eta(z)={(q^6z;q^4)_\infty\over(q^4z;q^4)_\infty},
\qquad
(z;p)_\infty=\prod_{n=0}^\infty(1-p^nz).
\]
In~\rf{gen} we have set $p_j=0$ $(\hbox{if $j\equiv N\bmod2$})$,
$=1$ (otherwise).

\begin{prop}\label{prop:Frel}
The Laurent coefficients of the following belong to $\Ker \rh'$.
\begin{eqnarray}
&&F_{\vep_1,\ldots,\vep_j,\vep_{j+1},\ldots,\vep_N}
(z_1,\ldots,z_{j+1},z_j,\ldots,z_N)\nonumber\\
&&\quad -\tilde{R}_{j, j+1}(z_{j+1}/z_j)
F_{\vep_1,\ldots,\vep_j,\vep_{j+1},\ldots,\vep_N}
(z_1,\ldots,z_j,z_{j+1},\ldots,z_N).
\label{eqn:Fcom}
\end{eqnarray}
\end{prop}

Proposition~\ref{prop:Frel} implies that the vectors
$\rh'_0\Bigl(F_{\vep_1\cdots\vep_N,m_1\cdots m_N}\Bigr)$ with
$m_1\le\cdots\le m_N$ span the space $\H$.

Let us introduce a second filtration in $\Vh_N^{'(r)}$.
\begin{eqnarray*}
&&\Vh^{'(r)}_N\supset \cdots \supset \Vh^{'(r)}_N[[m]]
\supset \Vh^{'(r)}_N[[m+1]] \supset \cdots,
\\
&&\Vh^{'(r)}_N[[m]]={\rm cl.}
{\rm span}_\C \{ F_{\vep_1\cdots\vep_N,m_1\cdots m_N}\in \Vh^{'(r)}_N
\mid \max(m_1,\ldots,m_N)\ge m\,\}
\end{eqnarray*}
where cl stands for the closure in $\Vh'$.
We denote by $\Vh^{(r)}_N$ the completion of
$\Vh^{'(r)}_N$ with respect to this filtration
\[
\Vh_N^{(r)}=
\lim_{\longleftarrow\atop m}\Vh^{'(r)}_N/\Vh^{'(r)}_N[[m]],
\]
and set $\Vh_N=\oplus_{r\in\Z}\Vh^{(r)}_N$, $\Vh=\oplus_{N\geq 0}\Vh_N$.

{}From Proposition~\ref{prop:Frel} it follows that
\[
\rh'\Bigl(\Vh^{'(r)}_N[[m]]\Bigr)\H_s=0
\]
if $m+s>0$. Therefore the maps
$\rh'$ and $\rh'_0$ extend further to $\Vh$:
\[
\rh~:~\Vh\longrightarrow\End_\C\H,
\quad
\rh_0~:~\Vh\longrightarrow \H.
\]

The vectors
$\rh'_0\Bigl(F_{\vep_1\cdots\vep_N,m_1\cdots m_N}\Bigr)$
are subject further to two kinds of relations:
the highest weight condition~\rf{Fhwt} and the fusion relation~\rf{Ffus}.

\begin{prop}
The kernel of $\rh_0~$ contains
\begin{equation}
F_{\vep_1\cdots\vep_N,m_1\cdots m_N}\quad \hbox{with $m_j>0$ for some $j$.}
\label{eqn:Fhwt}
\end{equation}
\end{prop}

\begin{prop}
The Laurent coefficients of the following belong to $\Ker \rh$.
\begin{eqnarray}
&& F_{\vep_1,\ldots,\vep_N}(z_1,\ldots,z_N)\Bigr|_{z_{j+1}=q^{-2}z_j}
\nonumber\\
&&\quad -(-q)^{N-j+(\vep_j-1)/2}\delta_{\vep_j+\vep_{j+1},0}
\prod_{i=1}^{j-1}(z_i-q^2z_j)\prod_{i=j+2}^{N}(q^{-2}z_j-q^2z_i)
\nonumber\\
&&\qquad \times
F_{\vep_1,\ldots,\vep_{j-1},\vep_{j+2},\ldots,\vep_N}
(z_1,\ldots,z_{j-1},z_{j+2},\ldots,z_N).
\label{eqn:Ffus}
\end{eqnarray}
\end{prop}

Define $\tN$ (respectively $\cN$) to be the closure of the span of
elements~\rf{Fcom} and~\rf{Fhwt} (respectively~\rf{Fcom}, \rf{Fhwt}
and~\rf{Ffus}) in $\Vh$.

\begin{prop}
The following map induced from $\rh_0$
is a $\uqc$-linear isomorphism
\[
\hV/\cN\simeq\H.
\]
\end{prop}

\subsection{Affine Hecke Algebra}
To distinguish different actions of $U'$ on the same module,
we shall sometimes write $U'_{k=l}$ for $l\in\C$ to mean the
quotient of $U'$ modulo $q^c=q^l$.
We will construct an action of $U'_{k=0}$ on $\hV$
in such a way that the subspace $\tN$ is invariant under the action.
The main tools are the affine Hecke algebras
and the $L$ operators acting on $\hV$.

Define $S\in\End_\C V\otimes V$ by
\begin{eqnarray*}
&&Sv_\vep\otimes v_\vep=-q^{-1}v_\vep\otimes v_\vep,\\
&&Sv_+\otimes v_-=(q-q^{-1})v_+\otimes v_--v_-\otimes v_+,\\
&&Sv_-\otimes v_+=-v_+\otimes v_-.
\end{eqnarray*}
The operators $S_{j,j+1}=(\pi_j\otimes\pi_{j+1})(S)\in\End_\C
V^{\otimes N}$ ($j=1,\ldots,N-1$) satisfy the Hecke algebra relations
\begin{eqnarray}
&&S_{j,j+1}-S_{j,j+1}^{-1}=q-q^{-1},
\label{eqn:Hecke1}\\
&&S_{j,j+1}S_{k,k+1}=S_{k,k+1}S_{j,j+1}
\qquad (|j-k|>1),
\label{eqn:Hecke2}\\
&&S_{j,j+1}S_{j+1,j+2}S_{j,j+1}=S_{j+1,j+2}S_{j,j+1}S_{j+1,j+2}.
\label{eqn:Hecke3}
\end{eqnarray}
The $R$-matrix $\tilde R(z)$~\rf{Rtilde} can be written as
\begin{equation}
\tilde R(z)={Sz-S^{-1}\over qz-q^{-1}}.
\label{eqn:RS}
\end{equation}

Define
\begin{equation}
G_{j,k}^{\pm 1}={q^{-1}z_j-qz_k\over z_j-z_k}(K_{j,k}-1)+q^{\pm 1}.
\label{eqn:G}
\end{equation}
Here $K_{j,k}$ signifies the exchange of variables $z_j$ and $z_k$.
The $G_{j,j+1}$ ($j=1,\ldots,N-1$) also satisfy
the relations~\rf{Hecke1}--\rf{Hecke3}.

\begin{prop}
  In terms of $S_{j,j+1}$ and $G_{j,j+1}$, the relation~\rf{Fcom} is
  written as \be
  (S_{j,j+1}-G_{j,j+1})F_{\vep_1,\ldots,\vep_N}(z_1,\ldots,z_N)=0.
  \en
\end{prop}

The action of the Hecke algebra given via $G_{j,j+1}$
can be extended to that of the affine Hecke algebra.
Set
\begin{eqnarray}
Y_j&=&(G_{j,j+1}^{-1}K_{j,j+1})\cdots(G_{j,N}^{-1}K_{j,N})
p^{\vartheta_j}(K_{1,j}G_{1,j})\cdots(K_{j-1,j}G_{j-1,j})
\nonumber\\
&=&G_{j,j+1}^{-1}\cdots G_{N-1,N}^{-1}Z
G_{1,2}\cdots G_{j-1,j},
\lb{GENY}
\end{eqnarray}
where $Z=K_{1,2}K_{1,3}\cdots K_{1,N}p^{\vartheta_1}$
and $p^{\vartheta_j}$ denotes the scale operator
\[
p^{\vartheta_j}f(z_1,\ldots,z_N)=f(z_1,\ldots,pz_j,\ldots,z_N).
\]

The operators $G_{j,j+1}$ ($j=1,\ldots, N-1$) and $Y_j$ ($j=1,\ldots,
N$) satisfy the relations for the affine Hecke algebra $\hat{H}_N$.
Namely we have, in addition to~\rf{Hecke1}--\rf{Hecke3},
\begin{eqnarray*}
&&Y_jY_k=Y_kY_j, \\
&&G_{j,j+1} Y_j G_{j,j+1}= Y_{j+1},\\
&&[G_{j,j+1}, Y_k]=0, \qquad (j, j+1\not=k).
\end{eqnarray*}

We note that the symmetric polynomials in $Y_j$ $(j=1,\ldots,N)$
belong to the center of $\hat{H}_N$.

For an operator $X\in\End_\C\C[z_1,z_1^{-1},\ldots,z_N,z_N^{-1}]$
we define $\hat X\in\End_\C\hV_N$ by
$$
\sum \left(\hat X F_{\e_1\cdots\e_N;m_1,\ldots,m_N}\right)
z_1^{-m_1}\cdots z_N^{-m_N}
=
\sum F_{\e_1\cdots\e_N;m_1,\ldots,m_N}
X\left(z_1^{-m_1}\cdots z_N^{-m_N}\right).
$$

\subsection{$L$-operators}

In order to define $U'_{k=0}$-action on $\hV$ we use the
$L$-operator formalism.  Some of the basic facts on $L$-operators are
given in Appendix~\ref{app:Drinf}.

Noting that the operators $\hY_j$ ($j=1,\ldots,N$) commute
with each other, we set
\begin{eqnarray}
\bL_0(x)&=&\bL_{01}(x;q^{N-1}\hY^{-1}_1)\cdots
\bL_{0N}(x;q^{N-1}\hY^{-1}_N),
\label{eqn:Lop}\\
\bL_{0j}(x;a)&=&\frac{xS_{0j}^{-1}-qa S_{0j}}{qa-x}P_{0j}.  \nonumber
\end{eqnarray}
Define $l^\pm_{jk}[n]\in\End_\C\hV$ as in~\rf{L} by expanding
$\bL_0(x)$ in $x^{\pm1}$.  The Yang-Baxter equation for
$\bR(x)$~\rf{Rbar} implies
\[
\bR_{00'}(x_0/x_{0'})\bL_0(x_0)\bL_{0'}(x_{0'})
=\bL_{0'}(x_{0'})\bL_0(x_0)\bR_{00'}(x_0/x_{0'}).
\]
Therefore, we get a representation of $\tilde U'$.
We denote it by
\[
\pi^{(N)}:\tilde U'\rightarrow \End_\C\hV.
\]

On $\uqc$, $\pi^{(N)}$ is equal to the representation induced by the
opposite coproduct. We have (cf.~\cite{CP})
\begin{eqnarray}
&&\pi^{(N)}(e_0)=\sum_{j=1}^Nq^{N-1}\hat Y^{-1}_j\pi_j(f_1)
\pi_{j+1}(t_1^{-1})\cdots\pi_N(t_1^{-1}),\label{eqn:newe}\\
&&\pi^{(N)}(f_0)=\sum_{j=1}^N q^{1-N}\hat Y_j\pi_1(t_1)
\cdots \pi_{j-1}(t_1)\pi_j(e_1),\label{eqn:newf}\\
&&\pi^{(N)}(t_0)=\pi_{1}(t_1^{-1})\cdots\pi_N(t_1^{-1}).
\label{eqn:newt}
\end{eqnarray}
Here we consider
\[
e_1=\pmatrix{0&1\cr0&0\cr},f_1=\pmatrix{0&0\cr1&0\cr},
t_1=\pmatrix{q&0\cr0&q^{-1}\cr},
\]
in $\End_\C\vaf$ and embed them by
$\pi_j:\End_\C\vaf\rightarrow\End_\C\hV$.

Because $Y^{\pm1}_j$ maps polynomials in $z$ to polynomials,
the linear span of~\rf{Fhwt} is invariant by $\pi^{(N)}$.
Because we have
\begin{eqnarray*}
&&\bL_{0j}(x;q^{N-1}\hY^{-1}_j)\bL_{0j+1}(x;q^{N-1}\hY^{-1}_{j+1})
\bigl(S_{j,j+1}-\hat G_{j,j+1}\bigr)\\
&&\hbox{\hskip1in}=\bigl(S_{j,j+1}-\hat G_{j,j+1}\bigr)
\bL_{0j}(x;q^{N-1}\hY^{-1}_{j+1})\bL_{0j+1}(x;q^{N-1}\hY^{-1}_j),
\end{eqnarray*}
the linear span of~\rf{Fcom} is also invariant.
Therefore we have
\begin{prop}
The $U'_{k=0}$-action $\pi^{(N)}$ on $\hV$ induces
a $U'_{k=0}$-action on $\hV/\tilde\cN$.
\end{prop}

In~\cite{USC} we have further shown that
\begin{prop}
The $U'_{k=0}$-action on $\hV/\tilde\cN$ induces
a $U'_{k=0}$-action on $\hV/\cN$.
\end{prop}


\setcounter{equation}{0}
\setcounter{prop}{0}
\def\omegat{\tilde{\omega}}

\section{Highest weight vectors} \label{sec:hwv}

\subsection{Vectors $\omegat_{\lambda,N}$}
In this subsection we fix $N\ge 1$ and work with the space $\hV/\tN$.
Following the ideas of~\cite{BPS},
we shall introduce a family of highest weight vectors
$\omegat_{\lambda,N}$ with respect to
the level-$0$ action, and find their Drinfeld polynomials.
In what follows we shall write $z$ for $(z_1,\ldots,z_N)$.

Set $\gamma(z)=\prod_{j<k}(z_j-q^2z_k)$, and define a generating series
in  $\hV/\tN$,
\begin{eqnarray}
\Omega(z)&=&\gamma(z)^{-1}F_{+\cdots +}(z)
\label{eqn:Om}\\
&=&\frac{\prod_{j=1}^Nz_j^{(j-N-p_j+p_N)/2}}{\prod_{j<k}\xi(z_k/z_j)}
\,v_+(z_1)\otimes\cdots\otimes v_+(z_N).
\end{eqnarray}
Here $\xi(z)=(q^2z;q^4)_\infty/(q^4z;q^4)_\infty$.

The following proposition is a direct consequence of
the properties of the series $F_{\e_1,\ldots,\e_n}(z)$.
\begin{prop}
The coefficients of $\Omega(z)$ are well-defined in $\hV/\tN$.
$\Omega(z)$ is symmetric in
$z_1,\ldots,z_N$ and does not comprise negative powers of them.
\end{prop}

The Macdonald polynomials $P_\lambda(z;p,t)$ with $l(\lambda)\le N$
constitute a basis of symmetric polynomials in $z_1,\ldots,z_N$ (see
Appendix~\ref{app:Mac}).  Hence any symmetric formal power series has
a unique expansion in terms of them.  Define vectors
$\omegat_{\lambda,N}\in \hV/\tN$ by
\begin{equation}\label{eqn:om}
\Omega(z)=
\sum_{l(\lambda)\le N}
P_\lambda(z;q^4,q^2)\omegat_{\lambda,N}.
\end{equation}

We will determine the Drinfeld polynomial for the $U'_{k=0}$-module
generated by $\omegat_{\lambda,N}$. Below we drop the index of
$\bL_0(x)$ when there is no chance of confusion.

\begin{prop} Let $\bL(x)$ be as in~\rf{Lop}. Then
\begin{eqnarray}
\bL(x)\Omega(z)&=&
\pmatrix{\hat{A}(x)\cdot\Omega(z)& * \cr
0 & 1\cdot\Omega(z) \cr},
\label{eqn:ABA}\\
\hat{A}(x)&=&
q^N\frac{\hat{\Delta}(q^{-2}x^{-1})}{\hat{\Delta}(x^{-1})},
\qquad
\hat{\Delta}(u)=\prod_{j=1}^N(1-q^N\hY^{-1}_j u).
\label{eqn:delhat}
\end{eqnarray}
\end{prop}

The action of $\hat{\Delta}(u)$ on $\omegat_{\lambda,N}$ is given
as follows.

\begin{prop}We have
\begin{eqnarray}
\hat{\Delta}(u)\omegat_{\lambda,N}&=&
\Delta_{\lambda,N}(u)\omegat_{\lambda,N},
\label{eqn:Dln}\\
\Delta_{\lambda,N}(u)&=&\prod_{j=1}^N(1-q^{-4d_k+1}u),
\qquad
d_k=\lambda_k+\frac{1}{2}(N-k).
\label{eqn:dk}
\end{eqnarray}
\end{prop}

\pf
Since the action of $\hY_j$ is defined
as $Y_j$ on the series $F_{\e_1,\ldots,\e_N}(z)$, we have
\begin{eqnarray*}
\hat{\Delta}(u)\Omega(z)
&=&
\prod_{j=1}^N(1-q^N \hY^{-1}_j u)
\,\Omega(z)
\\
&=&\gamma(z)^{-1}
\prod_{j=1}^N(1-q^N \hY^{-1}_j u)
{}~F_{+\cdots +}(z)
\\
&=&
\left(\gamma(z)^{-1}\,\prod_{j=1}^N(1-q^N {Y}^{-1}_j u)\,
{\gamma(z)}\right)\,\Omega(z).
\end{eqnarray*}
Under the conjugation by $\gamma(z)$,
the operators $G_{j,j+1}$ and $Y_j$ change to
\begin{eqnarray}
g_{j,j+1}^{\pm 1}&=&\gamma(z)^{-1}\circ G_{j,j+1}^{\pm 1}\circ \gamma(z)
= \frac{qz_j-q^{-1}z_{j+1}}{z_j-z_{j+1}}(1-K_{j,j+1})-q^{\mp 1},
\nonumber\\
y_j&=&
(-q^2)^{-N+1}
\gamma(z)^{-1}\circ Y_j \circ \gamma(z)
\nonumber\\
&=&g_{j,j+1}^{-1}\cdots g_{N-1,N}^{-1}\,
Z\, g_{1,2}\cdots g_{j-1,j}.
\label{eqn:yj}
\end{eqnarray}
Here we have used $p=q^4$.
Substituting these into~\rf{om} we obtain
\begin{eqnarray*}
&&\sum_{l(\lambda)\le N}P_\lambda(z;q^4,q^2)
\hat{\Delta}(u)\omegat_{\lambda,N}
\\
&&~=
\sum_{l(\lambda)\le N}
\left(\prod_{j=1}^N\bigl(1+(-1)^Nq^{2-N}y_j^{-1} u\bigr)
P_\lambda(z;q^4,q^2)\right) \omegat_{\lambda,N}.
\end{eqnarray*}
On the other hand, the formulas~\rf{DP} and~\rf{Mac0}
yield
\[
\prod_{j=1}^N(1+y_ju)P_\lambda(z;p,q^2)
=
\prod_{j=1}^N\Bigl(1-(-1)^N q^{N-2j+1}p^{\lambda_j}u\Bigr)\,
P_\lambda(z;p,q^2).
\]
The eigenvalues of an arbitrary symmetric function in $y_j$ are determined
from this formula.
\qed

\medskip

{}From~\rf{ABA}--\rf{Dln}
we see that $\omegat_{\lambda,N}$ is a highest weight vector
whose Drinfeld polynomial is given by~\rf{dk}.

\subsection{Bosonic expression}

Consider now the image of $\omegat_{\lambda,N}$ in $\H$ and call it
$\omega_{\lambda,N}$:
\[
\hV/\tN\rightarrow \hV/\cN\simeq \H,
\qquad
\omegat_{\lambda,N}\mapsto \omega_{\lambda,N}.
\]
In the next section we shall study the submodules of $\H$
generated by $\omega_{\lambda,N}$.

As noted in~\cite{BPS} in the case $q=1$, these vectors
can be expressed neatly in terms of bosons.
Recall the bosonic realization of the level-$1$ integrable modules
\[
V(\Lambda_i)=\C[b_{-1},b_{-2},\ldots]\otimes
\left(\bigoplus_{n\in\Z}\C e^{\Lambda_i+n\alpha}\right).
\]
Here $b_n$'s stand for the standard bosonic oscillators
satisfying $[b_m,b_n]=m\delta_{m+n,0}$ ($m,n\in\Z\backslash\{0\}$) .
The highest weight vector of $V(\Lambda_i)$ is given by
$\ket{i}=1\otimes e^{\Lambda_i}$.
In this language the `$+$'-component of the  vertex operator
reads as
\begin{eqnarray*}
\tilde{\Phi}^*_+(z)&=&
\exp\left(\sum_{n=1}^\infty \frac{b_{-n}}{n}q^{2n}z^n\right)
\exp\left(-\sum_{n=1}^\infty \frac{b_{n}}{n(1+q^{2n})}
z^{-n}\right)
\\
&&\otimes e^{\alpha/2}\left(-qz\right)^{(\partial+I)/2}(-q)^I.
\end{eqnarray*}
Here $I$ acts as $i$ times the identity on $V(\Lambda_i)$ ($i=0,1$).
For the details see~\cite{JM}.

\medskip\noindent{\sl Remark.}\quad
Following~\cite{Etin} we have rescaled the bosons $a_n$
in~\cite{JM} as
\[
a_n=q^{-n/2}\frac{[n]}{n}b_n, \qquad
a_{-n}=q^{n/2}\frac{[2n]}{n}b_{-n}
\qquad (n>0).
\]
\qed
\medskip

Comparing the generating series with the duality property of
the Macdonald polynomials~\rf{dual}, we find that
\begin{equation}\label{eqn:hwvB}
\omega_{\lambda,N}=(-q)^M
\tilde{P}_{\lambda'}(s;q^2,q^4)
\otimes e^{\Lambda_i+N\alpha/2},
\qquad M=\frac{N(N-1)}{2}+\frac{3}{2}\left[\frac{N}{2}\right],
\end{equation}
wherein the variables $s_n$ are substituted by the oscillators
\[
s_n=(-1)^{n-1}q^{2n}b_{-n},
\]
and $\left[x\right]$ denotes the integer-part of $x$.
{}From~\rf{hwvB} it is clear that the vectors
$\omega_{\lambda,N}\in\H$ have the homogeneous degree
\begin{equation}\label{eqn:deg}
{\rm deg}\,\omega_{\lambda,N}=|\lambda|+\left(\frac{N}{2}\right)^2-\frac{i}{4}
\end{equation}
where $|\lambda|=\lambda_1+\cdots+\lambda_N$ and $N\equiv i~(\mod 2)$
($i=0,1$).

\setcounter{equation}{0}
\setcounter{prop}{0}

\section{Decomposition} \label{sec:decomp}
In this section, we study how the level-$1$ integrable
modules $V(\Lambda_i)$ $(i=0,1)$ decompose
with respect to the new level-$0$ action.
We assume $|q|<1$ and that $q$ is transcendental over $\Q$.

\subsection{Submodules}

Consider the subspace generated by $\omega_{\lambda,N}$'s,
\begin{equation}\label{eqn:direct}
\sum_{\lambda,N} W'_{\lambda,N}\,\subset \H,
\qquad
W'_{\lambda,N}=U'_{k=0}\omega_{\lambda,N},
\end{equation}
where $N\ge 0$ and $\lambda$ runs over partitions with $l(\lambda)\le N$.

As the first step we show:
\begin{prop}\label{prop:direct}
The sum~\rf{direct} is direct.
\end{prop}

\begin{lem}
Let $I_1$ be the operator acting on $\hV_N$ by
\begin{equation}\label{eqn:int}
I_1=\sum_{j=1}^N q^{N-1}\hY_j^{-1}
-\frac{N}{2}(1+q^2).
\end{equation}
Then it is well-defined on $\H$ and commutes with the action of $U'_{k=0}$.
\end{lem}

\pf
By the construction of the action of $U'_{k=0}$,
symmetric polynomials in $\hY_j$
act on $\hV_N/(\tN\cap\hV_N)$, commuting with $U'_{k=0}$.
For~\rf{int} to be well-defined on $\H$,
we must show that it is compatible with the fusion relation.
Let us set
\begin{equation}\label{eqn:int1}
I^{(N)'}_1=\sum_{j=1}^N q^{N-1}\left(\hY^{(N)}_j\right)^{-1}
\end{equation}
where the $N$-dependence is exhibited explicitly.
For the action to be well-defined, the relations must hold:
\begin{eqnarray*}
&&\left(I^{(N)'}_1
F_{\vep_1\cdots\vep_N}(z_1,\ldots,z_N)\right)\Bigl|_{z_N=q^{-2}z_{N-1}}
\\
&&=
\left(I^{(N-2)'}_1+1+q^2\right)
\left(F_{\vep_1\cdots\vep_N}(z_1,\ldots,z_N)\Bigl|_{z_N=q^{-2}z_{N-1}}\right).
\end{eqnarray*}
This can be verified in almost the same way as in the proof
of Proposition~12~(4.8) in~\cite{USC}.
\qed

\medskip\noindent
{\sl Proof of Proposition~\ref{prop:direct}.}\quad
Since $I_1$ commutes with $U'$, it acts as a scalar on $W'_{\lambda,N}$.
{}From the formula~\rf{Dln}, its eigenvalue is given by
\[
\sum_{j=1}^N q^{-4d_j}-\frac{N}{2}(1+q^2).
\]
Note that $d_1>\cdots>d_N\ge 0$.
Since we assume $q$ to be transcendental,
these eigenvalues are all distinct
when $(\lambda,N)$ runs over all possible pairs.
The proposition follows from this fact.
\qed

\subsection{Irreducible modules}

Let $W_{\lambda,N}$ denote the irreducible module which has the same
Drinfeld polynomial~\rf{Dln} as $\omega_{\lambda,N}$.
At this stage we can only say that
\begin{equation}\label{eqn:estim}
\dim W'_{\lambda,N}\ge \dim W_{\lambda,N}.
\end{equation}
This holds because the irreducible module has the minimal dimension
among the finite dimensional
highest weight modules with the same highest weight.
As the next step, we study the structure of $W_{\lambda,N}$.

Given $\lambda$ with $l(\lambda)\leq N $,
let us re-parametrize it
by two sets of non-negative integers $\{ n_i\}$ and $\{k_i\}$ as follows:
\begin{equation}\label{eqn:kn}
\lambda=(\underbrace{k_1,\ldots,k_1}_{n_1}, \underbrace{k_2,\ldots,k_2}_{n_2},
\ldots, \underbrace{k_m,\ldots,k_m}_{n_m} ).
\end{equation}
Here $m$ is a positive integer,  $k_1>k_2>\cdots>k_m\geq 0$,
and $\sum_{i=1}^m n_i=N$.
Set
\begin{equation}\label{eqn:str}
a_j=q^{-4k_j-2N+2(n_1+\cdots+n_{j-1})+n_j+1}.
\end{equation}
Then the Drinfeld polynomial $\Delta_{\lambda,N}(u)$
in~\rf{Dln} is factorized as
\[
\Delta_{\lambda,N}(u)=
\Delta_{n_1}(u;a_1)\Delta_{n_2}(u;a_2)\cdots \Delta_{n_m}(u;a_m),
\]
with
\[
\Delta_{n}(u;a)=\prod_{i=1}^{n}(1-q^{n-2(i-1)}a u).
\]
To each polynomial $\Delta_{n_j}(u;a_j)$, one can associate the
$(n_j+1)$-dimensional evaluation module  $W_{n_j}(a_j)$ of $U'_{k=0}$
(see Appendix~\ref{app:Drinf} for the conventions).
The parametrization~\rf{kn}--\rf{str}
is so designed that for any
distinct pair $(n_j,a_j)$, $(n_k,a_k)$ the tensor product
$W_{n_j}(a_j)\otimes W_{n_k}(a_k)$ is irreducible.
This implies~\cite{ChPr91}

\begin{prop}
\begin{equation}\label{eqn:Wln}
W_{\lambda, N}\simeq W_{n_1}(a_1)\otimes \cdots\otimes W_{n_m}(a_m).
\end{equation}
\end{prop}

\subsection{Characters}
For $i=0,1$, define
\begin{equation}
{\cal V}^{(i)}\equiv
\bigoplus_{N\equiv i\ ({\rm mod}2) \atop
l(\lambda)\le N}W'_{\lambda, N}~\subset \H.
\end{equation}

The third step is to estimate its character
\begin{equation}
\ch_{{\cal V}^{(i)}}(q)=\tr_{{\cal V}^{(i)}}\left(q^{-d}\right)
=\sum_{n\ge 0}\Bigl(\dim {\cal V}^{(i)}_{-n}\Bigr) q^n
\end{equation}
where $d$ denotes the homogeneous grading operator on $\H$.
Given formal series $f(q)=\sum f_nq^n$ and $g(q)=\sum g_nq^n$,
we write $f(q)\gg g(q)$ to mean $f_n\ge g_n$ for all $n$.

\begin{lem}
\begin{equation}\label{eqn:est}
\ch_{{\cal V}^{(i)}}(q)\gg
\sum_{N\equiv i\ ({\rm mod}2)\atop l(\lambda)\le N}
q^{\frac{N^2}{4}-\frac{i}{4}+|\lambda|}\prod_{j=1}^m (n_j+1).
\end{equation}
\end{lem}
\pf
The degree of $\omega_{\lambda,N}$ is given by~\rf{deg}.
Since the action of $U'_{k=0}$ commutes with $d$,
this value of the degree occurs with
multiplicity $\dim W'_{\lambda,N}$.
{}From~\rf{estim} and~\rf{Wln} we obtain
\[
\dim W'_{\lambda,N}\ge \dim W_{\lambda,N}=\prod_{i=1}^m(n_i+1).
\]
The estimate~\rf{est} follows from this.
\qed
\medskip

It is known~\cite{BPS,ANOT,NaYa} that the right hand side
of~\rf{est} gives the (homogeneous) character for the level-$1$
module $V(\Lambda_i)$.
This leads us to the main result.
\begin{thm}
We have $W'_{\lambda,N}\simeq W_{\lambda,N}$.
The level-$1$ integrable highest weight module $V(\Lambda_i)$
is completely reducible under the action of $U'_{k=0}$:
\begin{equation}
V(\Lambda_i)\simeq
\bigoplus_{N\equiv i\, (\mod 2) \atop l(\lambda)\le N}
W_{\lambda, N}
\qquad (i=0,1).
\end{equation}
\end{thm}

\medskip\noindent{\sl Remark.}\quad
The action of the $U_q(\slt)$ is common to both level-$0$ and
level-$1$ structures.
Denoting by $h$ the standard generator of the Cartan subalgebra of
$\slt$, we have therefore
\[
\tr_{W'_{\lambda,N}}\left(q^dz^h\right)
=
q^{\frac{N^2}{4}-\frac{i}{4}+|\lambda|}\prod_{j=1}^m\chi_{n_j}(z)
\]
where
\[
\chi_n(z)=\frac{z^{n+1}-z^{-(n+1)}}{z-z^{-1}}
=\sum_{n^+,n^-\in \Z_{\geq 0}\atop n^++n^-=n}z^{n^+-n^-}
\]
is the character of the $(n+1)$-dimensional
irreducible representation of $\slt$.
The full character is then given by
\[
\ch_{{\cal V}^{(i)}}(q,z)=
\sum_{n^+,n^-\in\Z_{\geq 0}\atop n^++n^-\equiv i\ ({\rm mod} 2) }
z^{n^+-n^-}\frac{q^{\frac{(n^++n^-)^2}{4}-\frac{i}{4}}}
{(q;q)_{n^+}(q;q)_{n^-}}
\]
with $(q;q)_n=(1-q)\cdots(1-q^n)$.
This agrees with eq.(5.2) in~\cite{Mel}.
\qed

\appendix
\setcounter{equation}{0}
\setcounter{prop}{0}

\section{Drinfeld polynomials} \label{app:Drinf}

The aim of this appendix is to summarize our conventions concerning
various generators of $U'=U'\left(\slth\right)$ and Drinfeld
polynomials.

\subsection{Algebra $\tilde{U}'$}
Following~\cite{RS,DF}, define
$\tilde{U}'=U_q'\left(\glth\right)$
to be the associative algebra with unit,
generated by the symbols
$l^\pm_{ij}[\pm n]$ ($i,j=1,2$, $n\ge 0$) and an invertible
central element $q^{c/2}$.
In terms of the matrix generating series
\begin{equation}\label{eqn:L}
L^\pm(x)=\sum_{\pm n\ge 0}x^n
\pmatrix{ l^{\pm}_{11}[n] &  l^{\pm}_{12}[n] \cr
          l^{\pm}_{21}[n] &  l^{\pm}_{22}[n] \cr},
\end{equation}
the defining relations of $\tilde{U}'$ read as follows.
\begin{eqnarray*}
&&l^+_{ii}[0]l^-_{ii}[0]=1 \quad (i=1,2),
\qquad l^+_{21}[0]=l^-_{12}[0]=0,
\\
&& \bR_{0{0'}}(x_0/x_{0'})L^\pm_0(x_0)L^\pm_{0'}(x_{0'})
=L^\pm_{0'}(x_{0'})L^\pm_0(x_0)\bR_{0{0'}}(x_0/x_{0'}),
\\
&& \bR_{0{0'}}(q^c x_0/x_{0'})L^+_0(x_0)L^-_{0'}(x_{0'})
=L^-_{0'}(x_{0'})L^+_0(x_0)\bR_{0{0'}}(q^{-c}x_0/x_{0'}),
\end{eqnarray*}
where
\begin{equation}
\bR(x)=
\left(\matrix{
1  &         &          &    \cr
&\displ{\frac{(1-x)q}{1-q^2 x}}&\displ{\frac{(1-q^2)}{1-q^2x}}& \cr
&\displ{\frac{(1-q^2)x}{1-q^2x}}&\displ{\frac{(1-x)q}{1-q^2x}}& \cr
 &          &          & 1 \cr}\right)~.
\lb{Rbar}
\end{equation}
The matrices~\rf{L} can be uniquely decomposed~\cite{DF} as
\begin{eqnarray*}
&&L^\pm(x)=
\pmatrix{1 & f^\pm(x) \cr 0 & 1\cr}
\pmatrix{k^\pm_1(x) & 0 \cr 0 & k^\pm_2(x) \cr}
\pmatrix{1 & 0 \cr e^\pm(x) & 1 \cr},
\\
&&k^\pm_i(x)=\sum_{\pm n\ge 0}k^\pm_{i,n}x^n,\qquad
e^\pm(x)=\sum_{\pm n\ge 0}e^\pm_{n}x^n,\qquad
f^\pm(x)=\sum_{\pm n\ge 0}f^\pm_{n}x^n,
\end{eqnarray*}
with $e^+_0=f^-_0=0$.
Define $\tilde{x}_k^\pm$, $\tilde{\psi}_k,\tilde{\varphi}_{-k}$ by
\begin{eqnarray*}
&&k^+_1(x){k^+_2(x)}^{-1}=\sum_{k=0}^\infty \tilde{\varphi}_{-k}x^k,
\qquad
k^-_1(x){k^-_2(x)}^{-1}=\sum_{k=0}^\infty \tilde{\psi}_{k}x^{-k},
\\
&&
\frac{-1}{q-q^{-1}}\left(e^+(q^{c/2}x)-e^-(q^{-c/2}x)\right)
=
\sum_{n\in\Z} \tilde{x}_n^{+} x^{-n},
\\
&&
\frac{-1}{q-q^{-1}}\left(f^+(q^{-c/2}x)-f^-(q^{c/2}x)\right)
=
\sum_{n\in\Z} \tilde{x}_n^{-} x^{-n}.
\end{eqnarray*}

Let $x_k^\pm$ ($k\in\Z$), $\psi_k,\varphi_{-k}$ ($k\in\Z_{\ge 0}$),
$K=q^h$ and $\gamma$ be the Drinfeld generators of $U'$
in the notation of~\cite{JM}.
\begin{prop}[\cite{DF}]\label{prop:DF}
$U'$ is isomorphic to a subalgebra of $\tilde{U}'$ by the map
\[
\psi_k\mapsto \tilde{\psi}_k, \qquad
\varphi_k\mapsto \tilde{\varphi}_k,
\qquad
x^\pm_n\mapsto \tilde{x}^{\pm}_n,
\qquad
\gamma\mapsto q^c.
\]
\end{prop}


Define $l^\pm_{ij}(x)=\sum_{\pm n\ge 0}x^n l^{\pm}_{ij}[n]$. The
coproduct on $\tilde{U}'$ is defined by
\[
\Delta'\left(l^\pm_{ij}(x)\right)=
\sum_{k=1,2}l^\pm_{ik}(q^{\pm c_2/2}x)
\otimes l^\pm_{kj}(q^{\mp c_1/2}x)
\]
where $c_1=c\otimes 1$ and $c_2=1\otimes c$.
In terms of the Chevalley generators
\begin{equation}\label{eqn:Chev}
e_0=x_1^-t_0,\quad e_1=x_0^+,\quad
f_0=t_0^{-1}x_{-1}^+,\quad f_1=x_0^-,\quad
t_0=\gamma K^{-1},\quad t_1=K,
\end{equation}
the coproduct induced on $U'$ reads
\[
\Delta'(e_i)=e_i\otimes t_i+1\otimes e_i,
\qquad
\Delta'(f_i)=f_i\otimes 1+t_i^{-1}\otimes f_i,
\qquad
\Delta'(t_i)=t_i\otimes t_i.
\]
The identification of $e_0,f_0$ in~\rf{Chev} differs from
the one in~\cite{JM} by a factor of $\gamma$. Note also that
$\Delta'=\Delta^{\rm{ \scriptstyle op}}$.

\subsection{Drinfeld polynomials}
Let us recall from~\cite{ChPr91} some basic facts about the finite
dimensional representations of $U'$. On a finite-dimensional
$U'$-module $\gamma$ acts as $1$.

Let $W$ be a finite dimensional $U'$-module.
For a polynomial $P(u)\in \C[u]$ satisfying $P(0)=1$,
let $d_k^\pm\in\C$ be the coefficients in the expansion
\begin{eqnarray*}
q^{{\rm deg}P}\frac{P(q^{-2}u)}{P(u)}
&=& \sum_{k\ge 0}d_k^+ u^k \qquad (\hbox{ as } u\longrightarrow 0),
\\
&=& \sum_{k\le 0}d_k^- u^k \qquad (\hbox{ as } u\longrightarrow \infty).
\end{eqnarray*}
A vector $w\in W$
is said to be a highest weight vector
with the Drinfeld polynomial $P(u)$ if the following hold:
\begin{eqnarray}
&&x_k^+ w=0 \quad (k\in \Z),
\label{eqn:xk}\\
&&\psi_k w=d_k^+w \quad (k\ge 0),\qquad
\varphi_k w=d_k^-w \quad (k\le 0).
\label{eqn:phipsi}
\end{eqnarray}
Associating to $w$ its Drinfeld polynomial, one obtains
a bijective correspondence between the equivalence class of
irreducible finite dimensional modules
and the set of polynomials normalized as $P(0)=1$.

Suppose the module $W$ is defined through~\rf{L}.
If $w\in W$ has the property
\[
L^\pm(x)w=\pmatrix{ A^\pm(x)w & * \cr 0 & D^\pm(x)w\cr}
\]
with some functions $A^\pm(x)$, $D^\pm(x)$, then
the conditions~\rf{xk},~\rf{phipsi} are satisfied.
The Drinfeld polynomial $P(u)$ is determined by
\begin{eqnarray*}
q^{{\rm deg}P}\frac{P(q^{-2}u)}{P(u)}
&=& \frac{A^+(u^{-1})}{D^+(u^{-1})} \qquad (\hbox{ as } u\longrightarrow 0),
\\
&=& \frac{A^-(u^{-1})}{D^-(u^{-1})}
 \qquad (\hbox{ as } u\longrightarrow \infty).
\end{eqnarray*}

\medskip\noindent
{\sl Example.} \quad Let $W_n$ be the $(n+1)$-dimensional
irreducible module of $U_q(\slt)$. By $W_n(a)$ we mean the
evaluation module with the parameter $a\in \C\backslash\{0\}$
whose $U'$-module structure is given by
\[
e_0=af_1,\quad f_0=a^{-1}e_1,\quad t_0=t_1^{-1}.
\]
In our convention, the corresponding Drinfeld polynomial is
\[
P_n(u;a)=(1-q^nau)(1-q^{n-2}au)\cdots(1-q^{-n+2}au).
\]
\medskip


\setcounter{equation}{0}
\setcounter{prop}{0}

\section{Macdonald polynomials} \label{app:Mac}

\subsection{Macdonald polynomials}

We recall here some facts concerning the Macdonald polynomials
which are used in the text. For details see~\cite{Mac,Macbk}.
By a partition we mean a finite sequence of non-increasing integers
$\lambda=(\lambda_1,\ldots,\lambda_N)$,
$\lambda_1\geq\cdots\geq \lambda_N\geq 0$.
Its length is $l(\lambda)=\max\{j\mid \lambda_j>0\}$.

Let $\Lambda_N$ be the ring of symmetric polynomials in $N$
variables $z_1,z_2,\ldots,z_N$ with coefficients in $\Q(p,t)$, where
$p,t$ are indeterminates.  (In the usual notation the letter $q$ is
used in place of $p$.)  As before write $z$ for $(z_1,z_2,\ldots,z_N)$.
The Macdonald polynomials $P_\lambda(z;p,t)$ are a certain basis of
$\Lambda_N$ indexed by partitions $\lambda$ with $l(\lambda)\le N$.

The Macdonald polynomials are eigenfunctions of
the Macdonald operators defined as follows.
For each $n=1,\ldots,N$, let
\begin{equation}\label{eqn:Mop}
D_N^n(p,t)=t^{\frac{n(n-1)}{2}}\sum_{I}\Bigl(\prod_{i\in I\atop j\not\in I}
\frac{tz_i-z_j}{z_i-z_j}\Bigr)\prod_{k\in I}p^{\vartheta_k}.
\end{equation}
Here the summation is taken over all $n$-element subsets $I$ of
$\{1,2,\ldots,N\}$, and
$p^{\vartheta_j}$ denotes the shift operator
\[
 (p^{\vartheta_j}f)(z_1,\ldots,z_N)=f(z_1,\ldots,pz_j,\ldots,z_N).
\]
Setting $D^0_N(p,t)=1$, we consider the generating function
\begin{equation}\label{eqn:gMop}
D(u;p,t)=\sum_{n=0}^N D_N^n(p,t)\ u^{n}.
\end{equation}

\begin{prop}
\begin{equation}\label{eqn:DP}
D(u;p,t)P_{\lambda}(z;p,t)=\prod_{j=1}^N(1+t^{N-j}p^{\lambda_j}u)
P_{\lambda}(z;p,t).
\end{equation}
In particular
\begin{equation}\label{eqn:D1}
D_N^1(p,t) P_{\lambda}(z;p,t)=\Bigl(\sum_{j=1}^Nt^{N-j}p^{\lambda_j}\Bigr)
P_{\lambda}(z;p,t).
\end{equation}
\end{prop}

Let
$\Lambda=\displaystyle\lim_{\longleftarrow }\Lambda_N$ be
the projective limit given via the restriction
$\Lambda_{N+1}\rightarrow \Lambda_N$, $P\mapsto P\bigr|_{z_{N+1}=0}$.
For each fixed $\lambda$, $P_\lambda(z;p,t)$
gives rise to a well-defined element in $\Lambda$.
We denote it by $\tilde{P}_\lambda(s;p,t)\in\Lambda$
viewed as a polynomial in the power sum $s_n=\sum_jz_j^n$, $n\ge 1$.
Rewriting the duality property (\cite{Mac}, eq.(3.12)) we have
\begin{prop}Let
$\lambda'$ denote the partition conjugate to $\lambda$
(given by transposing the corresponding Young diagram). Then
\begin{equation}\label{eqn:dual}
\exp\left(\sum_{n\ge 1}(-1)^{n-1}\frac{s_n}{n}\sum_{j=1}^N w_j^n\right)
=\sum_{l(\lambda)\le N}
P_{\lambda}(w;p,t)\tilde{P}_{\lambda'}(s;t,p),
\end{equation}
where the sum ranges over all partitions of length at most $N$.
\end{prop}

\subsection{A lemma on Macdonald operators}

Consider the operators $y_j$, $r_{ij}=K_{ij}g_{ij}$ given in~\rf{yj}:
\begin{eqnarray*}
y_j&=&r_{-}(j,N)p^{\vartheta_j}r_{+}(1,j),
\\
r_{-}(i,j)&=&r_{i\,i+1}^{-1}r_{i\,i+2}^{-1}\cdots r_{ij}^{-1},
\\
r_{+}(i,j)&=&r_{i\,j}r_{i+1\,j}\cdots r_{j-1\,j},
\\
r_{ij}&=&b_{ji}(1-K_{ij})-q^{-1}K_{ij},
\qquad
r_{ij}^{-1}=b_{ij}(1-K_{ij})-qK_{ij},
\\
b_{ij}&=&-\frac{qz_i-q^{-1}z_j}{z_i-z_j},
\qquad
b_{ij}+b_{ji}=-q-q^{-1}.
\end{eqnarray*}
Here we regard $p,q$ as independent (we do not restrict to $p=q^4$).
In this subsection we prove the following

\begin{prop}\label{prop:Mop}
For any symmetric polynomial $P\in\Lambda_N$, we have
\begin{equation}\label{eqn:Mac0}
\prod_{j=1}^N\left(1+y_ju\right)P
=D\bigl((-q)^{-N+1}u;p,q^2\bigr)P
\end{equation}
where $D(u;p,t)$ is given in~\rf{gMop}.
\end{prop}

Let $\cL$ be a homogeneous polynomial in $y_j$ of degree $n$.
When applied to a symmetric polynomial $P\in \Lambda_N$,
we can move all the $K_{ij}$ to the right and rewrite $\cL P$ in the form
\begin{equation}\label{eqn:normal}
\tL P=\sum_{1\leq j_1<\cdots<j_n\leq N} f_{j_1,\ldots,j_n}(z)
p^{\vartheta_{j_1}}\cdots p^{\vartheta_{j_n}}\,P
\end{equation}
with some rational functions $f_{j_1,\ldots, j_n}(z)$.
We shall refer to $\tL$ as the normal form of $\cL$.
The following shows that it is uniquely determined from $\cL$.

\begin{lem}
  Suppose an operator $\tL$ of the form~\rf{normal} satisfies $\tL
  P=0$ for all $P\in\Lambda_N$, identically in $p$.  Then
  $f_{j_1,\ldots, j_n}(z)=0$ for all $j_1,\ldots,j_n$.
\end{lem}

\medskip\noindent{\it Proof.\quad} It is helpful to extend the
definition of $f_{j_1, \ldots ,j_n}(z)$ for all
$j_1,\ldots,j_n\in\{1,\ldots,N\}$ by symmetry in $j_1,\ldots,j_n$, and
as $0$ whenever two indices coincide.

Take $P=\sum_{\sigma\in S_N}
z^{\lambda_1}_{\sigma(1)}\cdots z^{\lambda_n}_{\sigma(n)}$
where $\lambda_j$'s are arbitrary positive integers.
By assumption $\tL P=0$. Picking the coefficient of
$p^{\lambda_1+\cdots+\lambda_n}$
we find
\begin{equation}\label{eqn:vander}
\sum_{j_1,\ldots,j_n=1}^N f_{j_1,\ldots,j_n}(z)
z^{\lambda_1}_{j_1}\cdots z^{\lambda_n}_{j_n}=0.
\end{equation}
Now if $\sum_{i=1}^N z_i^\lambda g_i(z)=0$ holds for any
$\lambda=1,2,\ldots$, then $g_i(z)=0$ for all $i$.  Since
$\lambda_j$'s are arbitrary,~\rf{vander} implies
$\sum_{j_2,\ldots,j_n} f_{j_1,\ldots,j_n}(z) z^{\lambda_2}_{j_2}\cdots
z^{\lambda_n}_{j_n}=0$ for all $j_1$.  Repeating this process we
obtain the Lemma.  \qed \medskip

{}From the Lemma, we must have
$\tL=\tL^\sigma$ for any $\sigma\in S_N$, where
\[
\tL^\sigma =\sum_{1\leq j_1<\cdots<j_n\leq N}
f_{j_1,\ldots,j_n}(z_{\sigma(1)},\ldots,z_{\sigma(N)})
p^{\vartheta_{\sigma(j_1)}+\cdots+\vartheta_{\sigma(j_n)}}.
\]
Therefore $\tL$ is determined once $f_{1,2,\ldots, n}(z)$ is known.

For each $n=1,\ldots,N$ define
\[
\cL_n=\sum_{1\leq j_1<\cdots<j_n \leq N} y_{j_1}\cdots y_{j_n}.
\]
Then Proposition~\ref{prop:Mop} is equivalent to the following.
\begin{prop}
  Let $I=\{1,2,\ldots,N\}$. The normal form of $\cL_n$ is given by
\begin{equation}\label{eqn:Mac}
\tL_n=
\sum_{|I|=n}\prod_{i\in I \atop j\not\in I}b_{ij}
\prod_{i\in I}p^{\vartheta_i}
=(-q)^{-n(N-1)}D_N^n(p,q^2)
\end{equation}
where $D_N^n(p,t)$ signifies the Macdonald operator~\rf{Mop}.
\end{prop}

\medskip\noindent{\it Proof. \quad}
We prove~\rf{Mac} by induction on $n$. The case
$n=0$ is trivially true.

{}From the remark above it suffices to compute the coefficient
of $p^{\vartheta_1+\cdots+\vartheta_n}$ in $\tL_n$.
If $j_1<\cdots<j_n$, then we have
\begin{eqnarray}
y_{j_1}\cdots y_{j_n}&=&
r_{-}(j_1, N)\,p^{\vartheta_{j_1}}\,r_{-}(j_2, N)\,
p^{\vartheta_{j_2}}\cdots r_{-}(j_n, N)\,p^{\vartheta_{j_n}}
\nonumber\\
&&\times
r_{+}(1,j_1)r_{+}(2,j_2)\cdots r_{+}(n,j_n).
\label{eqn:a1}
\end{eqnarray}
Suppose further $2\leq j_1<\cdots<j_n$, and
set $y'_j=r_{-}(j,N)p^{\vartheta_j}r_{+}(2,j)$.
We then find, using $r_{ij}P=-q^{-1}P$ for $P\in \Lambda_N$, that
\begin{equation}\label{eqn:a2}
y_{j_1}\cdots y_{j_n} P=
(-q)^{-n}y'_{j_1}\cdots y'_{j_n} P.
\qquad (P\in \Lambda_N)
\end{equation}

Let us write $\cL_n=y_1 \cL'_{n-1}+\cL''_n$, where
\[
\cL'_{n-1}=\sum_{2\leq j_2<\cdots<j_n\leq N} y_{j_2}\cdots y_{j_n},
\qquad
\cL''_n=\sum_{2\leq j_1<\cdots<j_n \leq N} y_{j_1}\cdots y_{j_n}.
\]
{}From~\rf{a1} it is clear that
the normal forms of
$\cL'_{n-1}$, $\cL''_n$ do not contain terms
$\prod_{i\in I}p^{\vartheta_i}$
with $1\in I$.
To calculate the coefficient of $p^{\vartheta_1+\cdots+\vartheta_n}$,
it is therefore enough to consider
$y_1\cL'_{n-1}$.
By~\rf{a2} we may apply the induction hypothesis to $\cL'_{n-1}$,
to obtain
\begin{eqnarray*}
y_1\cL'_{n-1} P&=&r_{-}(1,n)\cL'''P\qquad (P\in \Lambda_N),
\\
\cL'''&=&(-q)^{-n+1}r_{1\,n+1}^{-1}\cdots r_{1\,N}^{-1}\,p^{\vartheta_1}\,
\sum_{I'}\prod_{i\in I'\atop j\not\in I'}b_{ij} \prod_{i\in I'}p^{\vartheta_i}.
\end{eqnarray*}
Here the sum $\sum_{I'}$ ranges over the subsets $I'\subset \{2,\ldots,N\}$
with $n-1$ elements.
Let us bring the $K_{ij}$ contained in $r_{ij}$ to the right
and pick the coefficient of $p^{\vartheta_1+\cdots+\vartheta_n}$
in  the normal form of $r_{-}(1,n)\cL'''$.
The operator $r_{-}(1,n)$ comprises permutations of $z_1,\ldots,z_n$ alone,
so the coefficient of $p^{\vartheta_1+\cdots+\vartheta_n}$
comes from that of $\cL'''$.
It is easy to see that its normal form is
\[
\tL'''=
(-q)^{-n+1}\prod_{j=n+1}^N b_{1j}\prod_{2\leq i\leq n\atop n<j\leq N}b_{ij}
\times \prod_{i=1}^np^{\vartheta_i}~+\cdots
\]
with $\cdots$ denoting terms containing $p^{\vartheta_i}$ with $i>n$.
The first term is symmetric in the indices $\{1,\ldots,n\}$, so
$r_{1\,j}^{-1}$ ($2\leq j\leq n$) acts as $-q$.
Hence the coefficient of $p^{\vartheta_1+\cdots+\vartheta_n}$ in $\tL$
is given by
\[
(-q)^{n-1}(-q)^{-n+1}\prod_{1\leq i\leq n\atop n<j\leq N}b_{ij}
=\prod_{1\leq i\leq n\atop n<j\leq N}b_{ij}.
\]
This completes the proof.
\qed


\ifx\undefined\bysame
\newcommand{\bysame}{\leavevmode\hbox to3em{\hrulefill}\,}
\fi

\end{document}